\documentclass[12pt,final]{article}
\usepackage[margin=2.5cm,top=1.5cm]{geometry}
\usepackage[all,arc]{xy}
\usepackage{mathrsfs}
\usepackage{amsmath,amsthm,amssymb,enumerate}
\usepackage{color}

\newtheorem{definition}{Definition}[section]

\newtheorem{theorem}[definition]{Theorem}

\newtheorem{remark}[definition]{Remark}

\newcommand{\Romannum}[1]{\uppercase\expandafter{\romannumeral #1}}

\numberwithin{equation}{section}

\newcommand\keywordsname{Key words}
\newcommand\AMSname{AMS subject classifications}

\newenvironment{@abssec}[1]{%
     \if@twocolumn
       \section*{#1}%
     \else
       \vspace{.05in}\footnotesize
       \parindent .2in
         {\upshape\bfseries #1. }\ignorespaces
     \fi}
     {\if@twocolumn\else\par\vspace{.1in}\fi}

\begin{document}

\title{On the sum of reciprocal  generalized Fibonacci numbers
\footnote{P. Yuan's research is supported by the NSF of China (Grant No. 11271142) and
the Guangdong Provincial Natural Science Foundation(Grant No. S2012010009942).}}
\author{Pingzhi Yuan\footnote{{\it{Corresponding author:\;}}yuanpz@scnu.edu.cn.}, Zilong He \footnote{{\it{Email address:\;}}hzldew@qq.com.}, Junyi Zhou  }
\vskip.2cm
\date{{\small
School of Mathematical Sciences, South China Normal University,\\
Guangzhou, 510631, P.R. China\\
}}
\maketitle

\begin{abstract}
In this paper, we consider infinite sums derived from the reciprocals of the generalized Fibonacci numbers. We obtain some new and interesting identities for the generalized Fibonacci numbers.
\end{abstract}

\section{Introduction}
\hskip.6cm For any integer $n\ge0$, the famous Fibonacci numbers $F_n$ and Pell numbers are defined by the second-order linear recurrence sequences
$$F_{n+2}=F_{n+1}+F_n, \quad F_0=0, \,\,\, F_1=1$$
and
$$P_{n+2}=2P_{n+1}+P_n, \quad P_0=0, \,\,\, P_1=1.$$
There are many interesting results on the properties of these two sequences, see [1-9]. In 2009, Ohtsuka and Nakamura \cite{ON89} studied the properties of the Fibonacci numbers, and proved the following two interesting identities:
 $$\left\lfloor\left(\sum_{k=n}^\infty\frac{1}{F_k}\right)^{-1}\right\rfloor=\left\{
   \begin{array}{cc}
   F_{n-2}, &  \hskip.2cm   { \rm if}\,  n\, { \rm is\, even\, and} \, n\ge2; \\
    F_{n-2}-1,  &{\rm if}\,  n\, {\rm is\, odd\, and}\,  n\ge1.
   \end{array} \right.$$

$$\left\lfloor\left(\sum_{k=n}^\infty\frac{1}{F_k^2}\right)^{-1}\right\rfloor=\left\{
   \begin{array}{cc}
   F_{n-1}F_n-1, &  \hskip.2cm   { \rm if}\,  n\, { \rm is\, even\, and} \, n\ge2; \\
    F_{n-1}F_n,  &{\rm if}\,  n\, {\rm is\, odd\, and}\,  n\ge1,
   \end{array} \right.$$
where $\lfloor x\rfloor$ is the floor function, that is, it denotes the greatest integer less than or equal to $x$.  Recently,   Holliday and  Komatsu \cite{HK11} (Theorems 3 and 4), Xu and Wang \cite{XW13} proved the following interesting identities for the Pell numbers:
 $$\left\lfloor\left(\sum_{k=n}^\infty\frac{1}{P_k}\right)^{-1}\right\rfloor=\left\{
   \begin{array}{cc}
   P_{n-1}+ P_{n-2}, &  \hskip.2cm   { \rm if}\,  n\, { \rm is\, even\, and} \, n\ge2; \\
    P_{n-1}+P_{n-2}-1,  &{\rm if}\,  n\, {\rm is\, odd\, and}\,  n\ge1,
   \end{array} \right.$$
   $$\left\lfloor\left(\sum_{k=n}^\infty\frac{1}{P_k^2}\right)^{-1}\right\rfloor=\left\{
   \begin{array}{cc}
   2P_{n-1}+ P_{n}-1, &  \hskip.2cm   { \rm if}\,  n\, { \rm is\, even\, and} \, n\ge2; \\
    2P_{n-1}+P_{n},  &{\rm if}\,  n\, {\rm is\, odd\, and}\,  n\ge1,
   \end{array} \right.$$
   $$\left\lfloor\left(\sum_{k=n}^\infty\frac{1}{P_k^3}\right)^{-1}\right\rfloor=\left\{
   \begin{array}{cc}
   P_n^2P_{n-1}+3P_nP_{n-1}^2+\lfloor-\frac{61}{82}P_n-\frac{91}{82}P_{n-1}\rfloor, &  \hskip.2cm   { \rm if}\,  n\, { \rm is\, even\, and} \, n\ge2; \\
   P_n^2P_{n-1}+3P_nP_{n-1}^2+\lfloor\frac{61}{82}P_n+\frac{91}{82}P_{n-1}\rfloor,  &{\rm if}\,  n\, {\rm is\, odd\, and}\,  n\ge1,
   \end{array} \right.$$
 where providing $P_{-1}=P_1=1$. In \cite{ZW12} and \cite{XW13}, the authors asked whether there exists a computational formula for
  $$\left\lfloor\left(\sum_{k=n}^\infty\frac{1}{P_k^s}\right)^{-1}\right\rfloor,$$
 where $s\ge4$ is a positive integers.

  Let $p$ and $q$ be  integers such that $p^2+4q>0$. Define the generalized Fibonacci sequence $\{U_n(p, q)\}$, briefly $\{U_n\}$,
   as shown: for $n\ge2$
 $$U_n=pU_{n-1}+qU_{n-2},$$
 where $U_0=0, U_1=1$.
  The Binet formulae for $\{U_n\}$ is
 \begin{equation}\label{eq1}U_n=\frac{\alpha^n-\beta^n}{\alpha-\beta},\end{equation} 
 where $\alpha, \beta=(p\pm\sqrt{p^2+4q})/2$.

 The main purpose of this paper related to the computing problem of
 \begin{equation} U(s, n)=\left\lfloor\left(\sum_{k=n}^\infty\frac{1}{U_k^s}\right)^{-1}\right\rfloor\end{equation}
 for  $s=3$ and $q=-1$. For easy computation, we assume that $p=a$ is a positive integer and $q=-1$ throughout the paper. We have

\begin{theorem} Let $a\ge 3$ be a positive integer, and let $G_{n}$ be defined by the second-order linear  recurrence sequence $G_{n+2}=aG_{n+1}-G_n, \quad G_0=0, \,\,\, G_1=1$. Then for all $n\ge2$ we have
$$\left\lfloor\left(\sum_{k=n}^\infty\frac{1}{G_k^3}\right)^{-1}\right\rfloor=\left\{
   \begin{array}{cc}
   G_n^3-G_{n-1}^3-3\sum_{k=0}^{\lfloor\frac{n-4}{5}\rfloor}G_{n-3-5k}-2, &  \quad     a=3 \,\, \mbox{and} \,\, n\equiv3\pmod{5}; \\
   G_n^3-G_{n-1}^3-3\sum_{k=0}^{\lfloor\frac{n-4}{5}\rfloor}G_{n-3-5k}-1, &  \quad     \mbox{otherwise}.
   \end{array} \right.$$
\end{theorem}




 \section{Proof of the main result }

 In this section, we will prove our main result.
We  consider  the case that $\alpha\beta=1$ and $s=3$.

\begin{proof}
 From the Taylor series expansion of $(1-\varepsilon)^{-3}$ as $\varepsilon\to0$, we have
 $$(1-\varepsilon)^{-3}=1+\sum_{n=1}^\infty\frac{(n+1)(n+2)}{2}\varepsilon^n=1+3\varepsilon+6\varepsilon^2+O(\varepsilon^3).$$
Using (\ref{eq1}), we have
 $$\frac{1}{G_k^3}=\frac{(\alpha-\beta)^3}{\alpha^{3k}}\left(1-\frac{1}{\alpha^{2k}}\right)^{-3}$$
 $$=\frac{(\alpha-\beta)^3}{\alpha^{3k}}\frac{1}{\left(1-\frac{3}{\alpha^{2k}}+\frac{3}{\alpha^{4k}}+\frac{1}{\alpha^{6k}}\right)}$$
$$= \frac{(\alpha-\beta)^3}{\alpha^{3k}}\left[1+\frac{3}{\alpha^{2k}}+\frac{6}{\alpha^{4k}}+\frac{10\alpha^{4k}-15\alpha^{2k}+6}{\alpha^{4k}(\alpha^{2k}-1)^3}\right]$$
$$= (\alpha-\beta)^3\left[\frac{1}{\alpha^{3k}}+\frac{3}{\alpha^{5k}}+\frac{6}{\alpha^{7k}}+\frac{10\alpha^{4k}-15\alpha^{2k}+6}{\alpha^{7k}(\alpha^{2k}-1)^3}\right].$$

It is easy to check that
$$\frac{10}{\alpha^{9k}}<\frac{10\alpha^{4k}-15\alpha^{2k}+6}{\alpha^{7k}(\alpha^{2k}-1)^3}<\frac{11}{\alpha^{9k}}$$holds for $a\ge3$ and $k\ge2$.

Thus
 $$ \sum_{k=n}^\infty\frac{1}{G_k^3}=(\alpha-\beta)^3\left[\frac{1}{\alpha^{3n}}\cdot\frac{\alpha^3}{\alpha^3-1}+\frac{3}{\alpha^{5n}}\cdot\frac{\alpha^5}{\alpha^5-1}+
 \frac{6}{\alpha^{7n}}\cdot\frac{\alpha^7}{\alpha^7-1}+\sum_{k=n}^\infty\frac{10\alpha^{4k}-15\alpha^{2k}+6}{\alpha^{7k}(\alpha^{2k}-1)^3}\right]$$
 $$=\frac{(\alpha-\beta)^3\alpha^3}{\alpha^{3n}(\alpha^3-1)}\left[1+\frac{3}{\alpha^{2n}}\frac{\alpha^2(\alpha^3-1)}{\alpha^5-1}+\frac{6\alpha^4}{\alpha^{4n}}
 \frac{(\alpha^3-1)}{\alpha^7-1}+R_n\right],$$
 where
 $$R_n=\frac{(\alpha^3-1)\alpha^{3n}}{\alpha^3}\sum_{k=n}^\infty\frac{10\alpha^{4k}-15\alpha^{2k}+6}{\alpha^{7k}(\alpha^{2k}-1)^3}.$$
 Since $\sum_{k=n}^\infty\frac{1}{\alpha^{9k}}=\frac{\alpha^9}{\alpha^{9n}(\alpha^9-1)}$, we have
 $$\frac{10\alpha^6}{\alpha^{6n}(\alpha^6+\alpha^3+1)}<R_n<\frac{11\alpha^6}{\alpha^{6n}(\alpha^6+\alpha^3+1)}$$holds for $a\ge3$ and $k\ge2$.

Taking reciprocal, we get
 $$ \left(\sum_{k=n}^\infty\frac{1}{G_k^3}\right)^{-1}=
 \frac{(\alpha^3-1)\alpha^{3n}}{(\alpha-\beta)^3\alpha^3}\frac{1}{1+\frac{3}{\alpha^{2n}}\frac{\alpha^2(\alpha^3-1)}{\alpha^5-1}+\frac{6\alpha^4}{\alpha^{4n}}
 \frac{(\alpha^3-1)}{\alpha^7-1}+R_n}$$
 $$<\frac{(\alpha^3-1)\alpha^{3n}}{(\alpha-\beta)^3\alpha^3}\frac{1}{1+\frac{3}{\alpha^{2n}}\frac{\alpha^2(\alpha^3-1)}{\alpha^5-1}+\frac{6\alpha^4}{\alpha^{4n}}
 \frac{(\alpha^3-1)}{\alpha^7-1}+\frac{10\alpha^6}{\alpha^{6n}(\alpha^6+\alpha^3+1)}}$$
 $$<\frac{(\alpha^3-1)\alpha^{3n}}{(\alpha-\beta)^3\alpha^3}-\frac{3(\alpha^3-1)^2\alpha^{n}}{(\alpha-\beta)^3\alpha(\alpha^5-1)}+\delta_1,$$
 where
 $$\delta_1=-\frac{6\alpha(\alpha^3-1)^2}{\alpha^{n}(\alpha-\beta)^3(\alpha^7-1)}-\frac{10\alpha^3(\alpha^3-1)}{\alpha^{3n}(\alpha-\beta)^3(\alpha^6+\alpha^3+1)}
 +\frac{9\alpha(\alpha^3-1)^3}{\alpha^{n}(\alpha-\beta)^3(\alpha^5-1)^2}$$
 $$+\frac{36\alpha^3(\alpha^3-1)^3}{\alpha^{3n}(\alpha-\beta)^3(\alpha^5-1)(\alpha^7-1)}
 +\frac{36\alpha^5(\alpha^3-1)^3}{\alpha^{5n}(\alpha-\beta)^3(\alpha^7-1)^2}+\frac{60\alpha^5(\alpha^3-1)^2}{\alpha^{5n}(\alpha-\beta)^3(\alpha^5-1)(\alpha^6+\alpha^3+1)}$$
 $$+\frac{120\alpha^7(\alpha^3-1)^2}{\alpha^{7n}(\alpha-\beta)^3(\alpha^7-1)(\alpha^6+\alpha^3+1)}+\frac{100\alpha^9(\alpha^3-1)}{\alpha^{9n}(\alpha-\beta)^3(\alpha^6+\alpha^3+1)^2}$$
 since
 $$\frac{1}{1+\varepsilon}=1-\varepsilon+\varepsilon^2-\frac{\varepsilon^3}{1+\varepsilon}.$$
 An easy calculation shows that $\delta_1\le\frac{4}{\alpha^{n+3}}$ holds for $a\ge3$ and $k\ge2$. Therefore,
 $$ \left(\sum_{k=n}^\infty\frac{1}{G_k^3}\right)^{-1}<\frac{(\alpha^3-1)\alpha^{3n}}{(\alpha-\beta)^3\alpha^3}-\frac{3(\alpha^3-1)^2\alpha^{n}}{(\alpha-\beta)^3\alpha(\alpha^5-1)}+\delta_1$$
 $$\le\frac{\alpha^{3n}-\alpha^{3n-3}}{(\alpha-\beta)^3}-\frac{3(\alpha^3-1)^2\alpha^{n}}{(\alpha-\beta)^3\alpha(\alpha^5-1)}+\frac{4}{\alpha^{n+3}}$$
 $$=G_n^3-G_{n-1}^3-\frac{3\alpha^{n+2}}{(\alpha-\beta)\alpha(\alpha^5-1)}+\frac{3(\alpha-1)}{(\alpha-\beta)^3\alpha^n}-\frac{\alpha^3-1}{(\alpha-\beta)^3\alpha^{3n}}+\frac{4}{\alpha^{n+3}}$$
 $$<G_n^3-G_{n-1}^3-\frac{3\alpha^{n+2}}{(\alpha-\beta)\alpha(\alpha^5-1)}+\frac{3(\alpha-1)}{(\alpha-\beta)^3\alpha^n}+\frac{4}{\alpha^{n+3}}$$
 $$=G_n^3-G_{n-1}^3-\frac{3\alpha^{n+2}}{(\alpha-\beta)\alpha(\alpha^5-1)}+\lambda_1,$$
 where
 $$\lambda_1=\frac{3(\alpha-1)}{(\alpha-\beta)^3\alpha^n}+\frac{4}{\alpha^{n+3}}<0.1681$$for $a\ge3$ and $n\ge2$.

 Similarly,  we have
 $$ \left(\sum_{k=n}^\infty\frac{1}{G_k^3}\right)^{-1}=
 \frac{(\alpha^3-1)\alpha^{3n}}{(\alpha-\beta)^3\alpha^3}\frac{1}{1+\frac{3}{\alpha^{2n}}\frac{\alpha^2(\alpha^3-1)}{\alpha^5-1}+\frac{6\alpha^4}{\alpha^{4n}}
 \frac{(\alpha^3-1)}{\alpha^7-1}+R_n}$$
 $$>\frac{(\alpha^3-1)\alpha^{3n}}{(\alpha-\beta)^3\alpha^3}\frac{1}{1+\frac{3}{\alpha^{2n}}\frac{\alpha^2(\alpha^3-1)}{\alpha^5-1}+\frac{6\alpha^4}{\alpha^{4n}}
 \frac{(\alpha^3-1)}{\alpha^7-1}+\frac{11\alpha^6}{\alpha^{6n}(\alpha^6+\alpha^3+1)}}$$
 $$>G_n^3-G_{n-1}^3-\frac{3\alpha^{n+2}}{(\alpha-\beta)\alpha(\alpha^5-1)}+\lambda_2.$$
 Since
 $$\frac{1}{1+\varepsilon}=1-\varepsilon+\varepsilon^2-\varepsilon^3+\frac{\varepsilon^4}{1+\varepsilon},$$ and $\varepsilon=\frac{3}{\alpha^{2n}}\frac{\alpha^2(\alpha^3-1)}{\alpha^5-1}+\frac{6\alpha^4}{\alpha^{4n}}
 \frac{(\alpha^3-1)}{\alpha^7-1}+\frac{11\alpha^6}{\alpha^{6n}(\alpha^6+\alpha^3+1)}<0.3$ for $a\ge3$ and $n\ge2$, we have $\varepsilon^2-\varepsilon^3>0.7\varepsilon^2$, whence
 we can take
 $$\lambda_2=\frac{3(\alpha-1)}{(\alpha-\beta)^3\alpha^n}-\frac{\alpha^3-1}{(\alpha-\beta)^3\alpha^{3n}}-\frac{6\alpha(\alpha^3-1)^2}{\alpha^{n}(\alpha-\beta)^3(\alpha^7-1)}$$
 $$-\frac{11\alpha^3(\alpha^3-1)}{\alpha^{3n}(\alpha-\beta)^3(\alpha^6+\alpha^3+1)}+\frac{6.3\alpha(\alpha^3-1)^3}{\alpha^{n}(\alpha-\beta)^3(\alpha^5-1)^2}>0$$for $a\ge3$ and $n\ge2$.

 Consequently, we have shown that
 \begin{equation}\label{eq-1}G_n^3-G_{n-1}^3-\frac{3\alpha^{n+2}}{(\alpha-\beta)\alpha(\alpha^5-1)}+\lambda_2<\left(\sum_{k=n}^\infty\frac{1}{G_k^3}\right)^{-1}<G_n^3-G_{n-1}^3-
 \frac{3\alpha^{n+2}}{(\alpha-\beta)\alpha(\alpha^5-1)}+\lambda_1,\end{equation}
 where $0<\lambda_2<\lambda_1<0.1681$ for $a\ge3$ and $n\ge2$, and $\lambda_1<0.0053$ for $a\ge4$ and $n\ge3$.

 Now the calculations show that
 $$\frac{\alpha^{n+2}}{(\alpha-\beta)\alpha(\alpha^5-1)}=\left\{
   \begin{array}{cc}
   G_{n-3}+G_{n-8}+\cdots+G_7+G_2-\frac{\alpha^2+\alpha^3}{(\alpha-\beta)(\alpha^5-1)}, &  \quad     n\equiv0\pmod{5}; \\
   G_{n-3}+G_{n-8}+\cdots+G_8+G_3-\frac{\alpha^2+\alpha^3}{(\alpha-\beta)(\alpha^5-1)}, &  \quad     n\equiv1\pmod{5}; \\
   G_{n-3}+G_{n-8}+\cdots+G_9+G_4-\frac{\alpha+\alpha^4}{(\alpha-\beta)(\alpha^5-1)}, &  \quad     n\equiv2\pmod{5}; \\
   G_{n-3}+G_{n-8}+\cdots+G_{10}+G_5-\frac{1+\alpha^5}{(\alpha-\beta)(\alpha^5-1)}, &  \quad     n\equiv3\pmod{5}; \\
   G_{n-3}+G_{n-8}+\cdots+G_6+G_1-\frac{\alpha+\alpha^4}{(\alpha-\beta)(\alpha^5-1)}, &  \quad     n\equiv4\pmod{5}.
   \end{array} \right.$$

  The calculations also show that  $\frac{3(\alpha^2+\alpha^3)}{(\alpha-\beta)(\alpha^5-1)}>\lambda_1$ for $a\ge3$ and $n\ge2$; $\frac{3(\alpha+\alpha^4)}{(\alpha-\beta)(\alpha^5-1)}>\lambda_1$ for $a\ge3$ and $n\ge2$; and $\frac{3(1+\alpha^5)}{(\alpha-\beta)(\alpha^5-1)}>\lambda_1+1$ for $a=3$ and $n\ge3$;
   $0.87<\frac{3(1+\alpha^5)}{(\alpha-\beta)(\alpha^5-1)}<1$ for $a>3$ and $n\ge3$. Combining the calculations and (\ref{eq-1}),  we obtain
   $$\left\lfloor\left(\sum_{k=n}^\infty\frac{1}{G_k^3}\right)^{-1}\right\rfloor=\left\{
   \begin{array}{cc}
   G_n^3-G_{n-1}^3-3\sum_{k=0}^{\lfloor\frac{n-4}{5}\rfloor}G_{n-3-5k}-2, &  \quad     a=3 \,\, \mbox{and} \,\, n\equiv3\pmod{5}; \\
   G_n^3-G_{n-1}^3-3\sum_{k=0}^{\lfloor\frac{n-4}{5}\rfloor}G_{n-3-5k}-1, &  \quad     \mbox{otherwise}.
   \end{array} \right.$$

  Therefore we have prove Theorem 1.1.\end{proof}

\begin{remark} We can also compute the cases $s>3$ or $q=1$, however, the computations are much more complicated. So we stop here.\end{remark}

\end{document}